\documentclass[12pt,reqno]{amsart}
\usepackage{amssymb}

\usepackage{graphicx}
\usepackage{amscd}
\usepackage[all]{xy}

\newtheorem{theorem}{Theorem}
\newtheorem{corollary}{Corollary}
\newtheorem{definition}{Definition}
\newtheorem{example}{Example}

\newtheorem{remark}{Remark}
\numberwithin{equation}{section}
 \numberwithin{remark}{section}
\numberwithin{definition}{section} \numberwithin{lemma}{section}
\numberwithin{proposition}{section}
 \numberwithin{example}{section}
 \numberwithin{corollary}{section}
\numberwithin{theorem}{section}

\begin{document}

\title[]{The Whitehead exact sequence and the classification problem of homotopy types}

\author{Mahmoud Benkhalifa}
\email{makhalifa@uqu.edu.sa}
\address{Department of Mathematics, Faculty of Applied Sciences, Umm Al-Qura University, Mekka, Saudi Arabia}


\begin{abstract}
This paper defines an invariant associated to Whitehead's certain exact sequence of a simply connected CW-complex which is much more
elementary  - and less powerful - than the boundary
invariant of Baues. Nevertheless, in good cases, it classifies the homotopy types of CW-complexes.
\end{abstract}

\maketitle

\thispagestyle{empty}

\section{Introduction}
\medskip

Classification of spaces (for our purpose we restrict ourselves to simply connected CW-complexes) is a major task of algebraic homotopy. From the first fundamental invariants (homotopy and homology groups) to today's developements such as operads, a wide panel of algebraic objects try to determine CW-complexes and their morphisms. Rational homotopy builds an equivalence of categories between simply connected spaces without torsion and algebraic categories easy to define and to work with (Therefore our problem tackled below is solved in rational homotopy, see for example \cite{BauL}).

Such a nice situation is out of reach for CW-complexes with torsion.
In this paper we accept to limit ourselves to specific morphisms. But we obtain a very simple criterium
to detect topological morphisms by algebraic datas.

The starting point is the Hurewicz morphism which connects homotopy
to homology; if $X$ is a CW-complexe, we denote it as usual by: $h_*
: \pi_*(X) \rightarrow H_*(X,\mathbb{Z})$. Whitehead \cite{Whi}
inserted it into a long exact sequence:
$$\cdots \rightarrow
H_{n+1}(X,\mathbb{Z})\overset{b_{n+1}}{\longrightarrow }%
\Gamma^{X} _{n}\overset{}{\longrightarrow} \pi
_{n}(X)\overset{h_{n}}{\longrightarrow }%
H_{n}(X,\mathbb{Z})\rightarrow \cdots$$
From this he obtained then a
good invariant for $4$-dimensional CW-complexes.

\medskip

The program of Whitehead was to extend these results to higher
dimensions. Many years after, Baues \cite{Bau2} took afresh the
problem and developed an elaborated theory. He mimicked the
Postnikov sequence in a categorical and homological setting via the
towers of categories. The fundamental step is a recursive
construction of CW-complexes, starting from the Whitehead certain
exact sequence. Actually one defines a category  from the
$n-1$-skeleton and a sophisticated algebraic "boundary invariant"
which overlaps informations derived  from Whitehead's sequence. This
invariant is built using a homotopical construction, the "principal
reduction". We have to notice that this construction is
theoretically defined, but in general non effective (hardly
reachable by direct calculations on examples).  Nevertheless the
whole theory is very nice, gives a theoretical recursive tool for
determining CW complexes and their maps from Whitehead's exact
sequence. At last Baues was able to give a complete sets of
invariants for CW-complexes with cells in a short range of
dimensions; we can view this last case as the first significative
generalization of Whitehead's results on $4$-dimensional complexes.

\medskip

Our program is a compromise between Whitehead's one and the elaborated results of Baues. To begin with, we suppose given (simply-connected) CW-complexes $X$ and $Y$, and a commutative ladder of maps between their respective Whitehead certain exact sequences.

\begin{picture}(300,90)(25,10)
\put(65,80){$\cdots \rightarrow H_{n+1}(X,\Bbb
Z)\overset{b_{n+1}}{\longrightarrow }\Gamma^{X} _{n}\longrightarrow
\pi_{n}(X)\longrightarrow H_{n}(X,\Bbb
Z)\overset{b_{n}}{\longrightarrow }\cdots$}
\put(235,76){$\vector(0,-1){46}$} \put(186,76){$\vector(0,-1){46}$}
 \put(65,20){$\cdots \rightarrow
H_{n+1}(Y,\Bbb Z)\overset{b'_{n+1}}{\longrightarrow }\Gamma^{Y}
_{n}\longrightarrow \pi_{n}(Y)\longrightarrow H_{n}(Y,\Bbb
Z)\overset{b'_{n}}{\longrightarrow }\cdots$}
\put(110,76){$\vector(0,-1){46}$} \put(111,50){\scriptsize
$f_{n+1}$} \put(291,50){\scriptsize  $f_{n}$}
\put(290,76){$\vector(0,-1){46}$} \put(188,50){\scriptsize
$\gamma_{n}$} \put(236,50){\scriptsize $\Omega_{n}$}
\end{picture}

\noindent We add a collection of extensions  belonging to  a
homotopy invariant set (the set of the characteristic
$n$-extensions) to these data and suppose the above ladder of maps
is compatible with these extensions (we say that the ladder is a
strong map). In general these data are not sufficient to define a
topological map $\alpha : X \rightarrow Y$ from the algebraic map
$f_* : H_*(X,\mathbb{Z}) \rightarrow H_*(Y,\mathbb{Z})$ in the
ladder. We need a recursive condition about the compatibility of
$\alpha$ with Whitehead's $\Gamma$-groups, and call morphism,
denoted by $(f_*, \gamma_*)$, such a ladder of maps between
Whitehead's exact sequences. The main
theorem sounds as follows:\\

{\bf Main theorem}: Let $X$ and $Y$ be two simply connected CW-complexes. Any strong morphism $(f_*, \gamma_*)$ from
Whitehead's certain exact sequence of $X$ towards $Y$ 's one gives rise to a map $\alpha : X \rightarrow Y$ such that $H_*(\alpha)=f_*$.
\medskip

Analogous discussions and theorems take place in algebraic
categories such as differential graded Lie algebras or differential
graded free chain algebras (notice that we can define a homotopy
theory for these categories such that they work "similarly" as
CW-complexes). These are simpler cases and lead to more powerful
theorems \cite{Benk0},  \cite{Benk3}, \cite{Benk1},  \cite{B3},  \cite{B2}, \cite{Benk}, \cite{Ben1},  \cite{Ben2}, \cite{Ben3}, \cite{Ben4} and \cite{Ben7}

\smallskip

\noindent  The paper is organized as follows.

\noindent In section 2, we recall the basic definitions of
Whitehead's certain exact sequence and his theorem about
$4$-dimensional simply-connected CW-complexes and section 3, we
define the characteristic $n$-extensions. In section 4, we formulate
and prove the main theorem.

\section{The certain exact sequence of Whitehead}
\medskip

\subsection{The cellular complex and the Hurewicz morphism}
Let $X$ be a simply connected CW-complex defined by the collection
of its skeleta $(X_n)_{n\geq 0}$, where we can suppose
$X_0=X_1=\star$.

\noindent The long exact sequence of the pair $(X_n, X_{n-1})$ in
homotopy and in homology are connected by the Hurewicz morphism
$h_*$:
\[\xymatrix{
    \cdots\ar[r]^{\hspace{-6mm} i_{m,n}}&\pi_{m}(X_{n}) \ar[r]^{\hspace{-6mm} j_{m,n}} \ar[d]_{h_{m}}  & \pi_{m}(X_{n},X_{n-1})\ar[r]^{\beta_{m,n}}\ar[d]^{h_{m}}& \pi_{m-1}(X_{n-1}) \ar[d]^{h_{m-1}\mbox{\hspace{3.8cm}} (1)}\ar[r]&\cdots \\
     \cdots\ar[r]^{\hspace{-6mm} i^{H}_{m,n}}&H_{m}(X_{n},\Bbb Z) \ar[r]^{\hspace{-6mm} j^{H}_{m,n}} &
    H_{m}((X_{n},X_{n-1}),\Bbb Z)\ar[r]^{\beta^{H}_{m,n}}& H_{m-1}(X_{n-1},\Bbb
    Z)\ar[r]&\cdots
  }
\]

 \begin{remark}
 The following elementary  facts are well-known.
\begin{enumerate}
\label{r5}
\item The Hurewicz morphism $h_m : \pi_m(X_n, X_{n-1}) \rightarrow H_m((X_n, X_{n-1}),\mathbb{Z})$ is an isomorphism if $m\leq n$, non-trivial only if $m=n$.
\item $\pi_n(X_n, X_{n-1})$ is the free $\mathbb{Z}$-module generated by the n-cells of $X$.
\item $C_nX=\pi_n(X_n, X_{n-1})$ with the differential $d_n= j_n\circ\beta_n$, where $\beta_n= \beta_{n,n}$ and $j_n=j_{n,n}$, defines the cellular chain complex of $X$ (its homology is of course the singular homology $H_*(X)$). From now on, we omit to refer to $\mathbb{Z}$, it is understood that we deal only with integral homology. Moreover $\beta_n : C_nX \rightarrow \pi_{n-1}(X^{n-1})$ represents by adjunction the attaching map for the $n$-cells $\vee S^n \rightarrow X^{n-1}$.
\end{enumerate}
 \end{remark}

\subsection{The definition of Whitehead's certain exact sequence}
Now Whitehead \cite{Whi} inserted the Hurewicz morphism in a long exact sequence connecting homology and homotopy. First he defined the following group
\begin{equation}\label{2}
\Gamma^X_n=\mbox{Im }( i_n : \pi_n(X_{n-1}) \rightarrow
\pi_n(X_n))=\mbox{ker }j_n , \forall n\geq 2.
\end{equation}
We notice that $\beta_{n+1}\circ d_{n+1}=0$ and so $\beta_{n+1}:
\pi_{n+1}(X_{n+1}, X_n)\rightarrow  \pi_n(X_n)$ factors through the
quotient: $b_{n+1}:  H_{n+1}(X)\rightarrow  \Gamma^X_n$.

\noindent With this map, Whitehead  \cite{Whi} defined the following
sequence:
\begin{equation}
\label{3} \cdots \rightarrow
H_{n+1}(X,\mathbb{Z})\overset{b_{n+1}}{\longrightarrow }%
\Gamma^{X} _{n}\overset{}{\longrightarrow} \pi
_{n}(X)\overset{h_{n}}{\longrightarrow }%
H_{n}(X,\mathbb{Z})\rightarrow \cdots
\end{equation}
 and proved the following.
\begin{theorem}
The above sequence is a natural exact sequence, called the certain exact sequence.
\end{theorem}

\noindent Notation. We shall denote the sequence (\ref{3}) by
$\mbox{WES }(X)$.

 This sequence improves the informations provided by both
homology and homotopy groups. Whitehead was led to the very natural
question: for which class of CW-complexes and maps does the certain
exact sequence define a complete invariant? In other words, given
the following commutative diagram of group maps:

\begin{picture}(300,90)(25,10)
\put(65,80){$\cdots \rightarrow H_{n+1}(X,\Bbb
Z)\overset{b_{n+1}}{\longrightarrow }\Gamma^{X} _{n}\longrightarrow
\pi_{n}(X)\longrightarrow H_{n}(X,\Bbb
Z)\overset{b_{n}}{\longrightarrow }\cdots$}
\put(235,76){$\vector(0,-1){46}$} \put(186,76){$\vector(0,-1){46}$}
 \put(65,20){$\cdots \rightarrow
H_{n+1}(Y,\Bbb Z)\overset{b'_{n+1}}{\longrightarrow }\Gamma^{Y}
_{n}\longrightarrow \pi_{n}(Y)\longrightarrow H_{n}(Y,\Bbb
Z)\overset{b'_{n}}{\longrightarrow }\cdots$}
\put(110,76){$\vector(0,-1){46}$} \put(111,50){\scriptsize
$f_{n+1}$} \put(291,50){\scriptsize  $f_{n}$}
\put(290,76){$\vector(0,-1){46}$} \put(188,50){\scriptsize
$\gamma_{n}$} \put(236,50){\scriptsize $\Omega_{n}$} \put(428,50){
$(2)$}
\end{picture}

\noindent what can we say about the existence of a cellular map
$\alpha : X\rightarrow Y$ with $H_*(\alpha)=f_*$?

\noindent The question has no answer in general. Whitehead
\cite{Whi}  gave a complete answer in the case of $4$-dimensional
simply connected CW-complexes. We recall it right now. Baues
  gave a more general and sophisticated answer; it needs
 long definitions and  new formulations (see \cite{Bau2} for details).

\subsection{$4$-dimensional CW-complexes}

We need first to define an algebraic functor which represents quadratic maps.

\noindent A function $f :A \rightarrow B$ between abelian groups is
called a quadratic map if $f(-a) = a$ and if the function $ A \times
A \rightarrow B$, defined by $(a, b) \mapsto f(a+b)-f(a)-f(b)$ is a
bilinear map. The following assertion is both a definition and the
proposition justifying it:

\begin{definition}
For every abelian group $A$ there exists a universal quadratic map
$$ \gamma : A\rightarrow \Gamma(A)$$
such that  every quadratic map $f :A \rightarrow B$ uniquely factorizes

\begin{picture}(300,90)(-35,30)
\put(120,98){$\vector(0,-1){42}$} \put(124,98){$\vector(3,-2){63}$}
 \put(156,78){ \scriptsize $f$}
\put(110,102){$A$} \put(110,45){$\Gamma(A)\vector(1,0){50} B$}
\put(107,78){\scriptsize $\gamma$}
\end{picture}

\noindent For any morphism $\phi : A\rightarrow A'$, $\Gamma (\phi)
: \Gamma(A) \rightarrow \Gamma(A')$ is defined and $\Gamma$ is a
well-defined functor, called Whitehead's quadratic functor.
\end{definition}

This  functor has the following properties:

1) If $\eta : S^3 \rightarrow S^2$ is the Hopf map, the induced map $\eta^* : \pi_2(X) \rightarrow  \pi_3(X) $ is quadratic;

2) $\Gamma_3^X = \Gamma(\pi_2(X))$;

3) $\Gamma_{n+1}^X = \pi_n(X)
\otimes\mathbb{Z}/2\,\,\,\,\,,\,\,\,\,n\geq 3$.

\noindent Then we can formulate Whitehead's theorem on
$4$-dimensional CW-complexes:
\begin{theorem}
\label{t1} Let $X$ and $Y$ be two simply connected $4$-dimensional
CW-complexes. We suppose there exists a commutative ladder of group
maps from WES($X$) towards WES($Y$) (notice that the Hurewicz map is
an isomorphism in degree $2$, so we can shorten the exact
sequences).

\begin{picture}(300,90)(25,10)
\put(65,80){$H_{4}(X^4,\Bbb Z)\overset{b_{4}}{\longrightarrow
}\Gamma(H_{2}(X^4,\Bbb Z))\longrightarrow
\pi_{3}(X^4)\twoheadrightarrow H_{3}(X^4,\Bbb Z)$}
\put(245,76){$\vector(0,-1){46}$} \put(166,76){$\vector(0,-1){46}$}
 \put(65,18){$H_{4}(Y^4,\Bbb
Z)\overset{b_{4}}{\longrightarrow }\Gamma(H_{2}(Y^4,\Bbb
Z))\longrightarrow \pi_{3}(Y^4)\twoheadrightarrow H_{3}(Y^4,\Bbb
Z)$} \put(88,76){$\vector(0,-1){46}$} \put(89,50){\scriptsize
$f_{4}$} \put(303,50){\scriptsize  $f_{3}$}
\put(302,76){$\vector(0,-1){46}$} \put(169,50){\scriptsize
$\gamma_4$} \put(247,50){\scriptsize $\Omega_{3}$}
\end{picture}

\noindent If $\gamma_4=\Gamma(f_2)$, there exists a cellular map
$\alpha : X\rightarrow Y$ with $H_n(\alpha)=f_n$, $n=2, 3, 4$.
\end{theorem}

\noindent Such a simple theorem is not valid for CW-complexes of
higher dimensions. Nevertheless we remark that we enforced  that an
invariant  for $4$-skeleton (namely the $\Gamma$ group) can be
calculated with invariants of the $3$-skeleton (namely homology in
degree $2$).

The following sections presents an elementary approach of the
problem, easy to calculate with but less powerful - than the
boundary invariant of Baues. Nevertheless, in good cases, it
classifies the homotopy types of CW-complexes.

\section{The characteristic $n$-extensions}
In order to give a partial generalization to theorem \ref{t1}, we
shall see that we add two ingredients to our receipts. First some
analogue to the condition  $\gamma_4=\Gamma(f_2)$. Second we add a
homotopy invariant set   to Whitehead's certain exact sequence,
called the set of the characteristic $n$-extensions,  which expresse
a compatibility condition for morphisms. The purpose of this section
is to define this set.
\subsection{Splitting homotopy groups}
Our task is to recursively define maps between spaces from morphisms
of their Whitehead certain exact sequences.

\noindent Consider now the morphism $j_n : \pi_n(X)\rightarrow C_nX$
extracted from diagram (1). It gives rise to the short exact
sequence
\begin{equation}\label{12}
\Gamma^X_n \rightarrowtail  \pi_n(X^n) \twoheadrightarrow  \mbox{ker
}\beta_n=\mbox{Im }j_n.
\end{equation}
 As $C_nX$ is a free abelian group, $\mbox{ker }\beta_n\subset C_nX$ is also free and the later short exact sequence splits.
 So we can choose a splitting
\begin{equation}\label{13}
\mu_n : \pi_n(X^n) \stackrel{\cong}{\rightarrow}\Gamma^X_n \oplus
\mbox{ker }\beta_n.
\end{equation}
\begin{remark}\label{r2} Let us denote by $\mu^1_n : \pi_n(X^n)\rightarrow
\Gamma^X_n$ the composition of $\mu_n$ with the projection onto the
first factor while the composition with the second projection is
$j_n$ (writing this in a matrix setting: $\mu_n= \left(
\begin{array}{cc}
\mu^1_n \\
j_n \\
\end{array}
\right)
$). If $\mu^{-1}_n$ is the inverse of $\mu_n$, it is clear that $\mu^{-1}_n|_{\Gamma^X_n}$ is the inclusion $\Gamma^X_n \subset \pi_n(X^n) $. If we denote by $\sigma_n=\mu^{-1}_n|_{\mbox{ker }\beta_n}$ the section of $j_n$ defined by $\mu_n$, we have the identification
\begin{equation}\label{14}
\mu^1_n=id_{\pi_n(X^n)}-\sigma_n\circ j_n.
\end{equation}

\end{remark}
The sequence (\ref{12}) is natural  and induces a commutative
diagram for any map $\alpha : X\rightarrow Y$

\begin{picture}(300,90)(-30,30)
\put(90,102){$\Gamma^{X}_{n}\rightarrowtail
\pi_{n}(X^{n})\twoheadrightarrow \ker \beta_{n}$}
\put(90,48){$\Gamma^{Y}_{n}\rightarrowtail
\pi_{n}(Y^{n})\twoheadrightarrow \ker \beta'_{n}$}
 \put(94,98){$\vector(0,-1){39}$}
 \put(130,98){$\vector(0,-1){39}$}
 \put(132,78){\scriptsize $\pi_{n}(\alpha_{n})$}
\put(187,78){\scriptsize $C_{n}\alpha_{| \ker \beta_{n}}$}
\put(77,78){\scriptsize $\gamma^{\alpha_{n}}_{n}$}
\put(185,98){$\vector(0,-1){39}$} \put(378,78){$(3)$}
\end{picture}

\noindent where $\alpha_n$ is induced from $\alpha$ by restriction
to the $n$-skeleton, $\gamma_n^{\alpha_n}$ is the restriction of
$\pi_n(\alpha_n)$ to $\Gamma^X_n \subset \pi_n(X^n)$.
  Using the splitting $\mu_n$, we can form the following (non-commutative!)
  diagram:

\begin{picture}(300,90)(-30,30)
\put(90,102){$\pi_{n}(X^{n})\overset{\mu_{n}}{\longrightarrow}\Gamma^{X}_{n}\oplus
\ker \beta_{n}$}
\put(90,48){$\pi_{n}(Y^{n})\overset{\mu'_{n}}{\longrightarrow}\Gamma^{Y}_{n}\oplus
\ker \beta'_{n}$}
 \put(104,98){$\vector(0,-1){39}$}
\put(178,78){\scriptsize $\gamma^{\alpha_{n}}_{n}\oplus
C_{n}\alpha_{| \ker \beta_{n}}$} \put(73,78){\scriptsize
$\pi_{n}(\alpha_{n})$} \put(176,98){$\vector(0,-1){39}$}
\put(378,78){(4)}
\end{picture}

\noindent   Let us examine $(\gamma_n^{\alpha_n} \oplus
C_n\alpha|_{\mbox{ker }\beta_n}) \circ \mu_n - \mu'_n\circ
\pi_n(\alpha_n): \pi_n(X^n) \rightarrow \Gamma^Y_n \oplus \mbox{ker
}\beta'_n$. By the above remark, the second summand is
$C_n\alpha|_{\mbox{ker }\beta_n}\circ j_n - j'_n \circ
\pi_n(\alpha_n)$. By diagram (3) it is zero. Therefore:
$$\mbox{Im }[(\gamma^{\alpha_n} \oplus  C_n\alpha|_{\mbox{ker }\beta_n}) \circ \mu_n - \mu'_n\circ
\pi_n(\alpha_n): \pi_n(X^n) \rightarrow \Gamma^Y_n \oplus \mbox{ker }\beta'_n] \subset \Gamma^Y_n.$$

Moreover, using the decomposition of $\mu_n$, formula (\ref{14}) and
the commutative diagram (3), we have:
\begin{equation}\label{17}
(\gamma_n^{\alpha_n} \oplus  C_n\alpha|_{\mbox{ker }\beta_n})\circ
\mu_n - \mu'_n\circ \pi_n(\alpha_n) = \pi_n(\alpha_n) \circ \sigma_n
\circ j_n - \sigma'_n \circ j'_n \circ \pi_n(\alpha_n).
\end{equation}
\subsection{Splitting the cellular complex. The characteristic $n$-extensions}

Consider the differential of the cellular complex: $d_{n+1} :
C_{n+1}X \rightarrow C_nX$; $\mbox{Im }d_{n+1} \subset C_nX$ is a
free abelian group. We can choose a splitting:
\begin{equation}\label{18}
t_{n+1} :  \mbox{Im }d_{n+1} \oplus \mbox{ker }d_{n+1}
\stackrel{\cong}{\rightarrow} C_{n+1}X
\end{equation}
whose restriction  to $\ker d_{n+1}$ is the inclusion.

\noindent  We can now rewrite the morphism $\beta_{n+1} : C_{n+1}X
\rightarrow \pi_n(X^n)$ using the above respective splittings of
source (\ref{18}) and target (\ref{13}) groups:
\begin{equation}\label{19}
  \mbox{Im }d_{n+1} \oplus \mbox{ker }d_{n+1} \stackrel{t_{n+1}}{\rightarrow} C_{n+1}X
  \stackrel{\beta_{n+1}} {\rightarrow} \pi_n(X^n) \stackrel{\mu_n}{\rightarrow}  \Gamma^X_n \oplus \mbox{ker
  }\beta_n
\end{equation}
and we write down the composition as a matrix:
\begin{equation}\label{1111}
  \left(
\begin{array}{cc}
\phi_n & \theta_n \\
\psi_n & \eta_n \\
\end{array}
\right)
 \end{equation}
 First, using formulas (\ref{14}) and the argument (3) in
remark \ref{r5} we get:
\begin{eqnarray}
\label{20}
\phi_n =
 \mu^1_n  \circ \beta_{n+1}\circ t_{n+1}|_{\mbox{Im }d_{n+1}}&=&(id_{\pi_n(X^n)}-\sigma_n\circ j_n) \circ \beta_{n+1}\circ t_{n+1}|_{\mbox{Im }d_{n+1}}\nonumber\\
&=& (\beta_{n+1} - \sigma_n \circ d_{n+1})\circ t_{n+1}|_{\mbox{Im
}d_{n+1}}.
\end{eqnarray}
Second, using the same argument and the remark \ref{r2} we get:
$$\psi_n = d_{n+1}\circ t_{n+1}|_{\mbox{Im }d_{n+1}}.$$
Finally, by definition of $b_{n+1}$ given in (\ref{2}):
$$\theta_n=b_{n+1}\circ pr_{n+1},$$
(where $pr_{n+1} : \mbox{ker
}d_{n+1} \twoheadrightarrow H_{n+1}(X)$ is the projection) and
$\eta_n=0.$

\noindent  Among the four components of the matrix $ \left(
\begin{array}{cc}
\phi_n & \theta_n \\
\psi_n & \eta_n \\
\end{array}
\right) $, only $\phi_n$ reflects data non-directly readable in the
Whitehead exacts sequence. $\phi_n$ depends on our two splittings
$t_{n+1}$ and $\mu_n$. Define now $\widetilde{\phi}_n$ by composing
$\phi_n$ with the projection $\Gamma^X_n\rightarrow \mbox{Coker
}b_{n+1}$. Then remark that:
\begin{equation}\label{21}
\mbox{Im }d_{n+1}  \stackrel{\kappa_n}{\rightarrowtail}\mbox{ker
}d_n \rightarrow H_n(X)
\end{equation}
is a free resolution of $H_n(X)$ so $\widetilde{\phi}_n$  defines
the extension class:
 \begin{equation}\label{22}
 [\widetilde{\phi}_n] \in \mbox{Ext }^1_{\mathbb{Z}}(H_n(X), \mbox{Coker }b_{n+1}).
 \end{equation}
\begin{definition}
\label{d1} The class $[\widetilde{\phi}_n] $ is called a
characteristic $n$-extension of the CW-complex $X$.
\end{definition}
\begin{remark}
\label{r7}It is important to  notice the following fact:

\noindent Let $\mathrm{Ext }(\ker b_n, \mathrm{Coker
    }\,b_{n+1})$ be the abelian group of the extensions classes of  $\mathrm{Coker
    }\,b_{n+1}$ by $\ker b_n$.

    \noindent It is well-known that $\mathrm{Ext }(\ker b_n, \mathrm{Coker
    }\,b_{n+1})$ and $\mbox{Ext }^1_{\mathbb{Z}}(\ker b_n, \mathrm{Coker
    }b_{n+1})$ are isomorphic. Now if we consider the surjection:
\begin{equation}\label{M1}
    \theta_{n}:\mathrm{Ext }^1_{\Bbb Z}(H_n(X), \mathrm{Coker
    }\,b_{n+1})\twoheadrightarrow \mathrm{Ext }(\ker b_n, \mathrm{Coker }\,b_{n+1}).
\end{equation}
induced by the inclusion $\ker b_{n}\subset H_{n}(X)$ and the above
isomorphism,    we can say that  any characteristic $n$-extension of
 $X$ satisfies:
\begin{equation}\label{M2}
\theta_{n}([\widetilde{\phi}_n] )=[\pi_{n}(X)]
\end{equation}
where $[\pi_{n}(X)]$ is the class represented by the short exact
sequence $\mathrm{Coker }\,b_{n+1}\rightarrowtail
\pi_{n}(X)\twoheadrightarrow \ker b_n(X)$  extracted for the
Whitehead exact sequence (\ref{2}). Therefore we can say that
$S_{n}(X)=(\theta_{n})^{-1}([\pi_{n}(X)])$ is the set of all  the
characteristic $n$-extensions of the CW-complex $X$.

\noindent The set $S_{n}(X)$ is an invariant of homotopy. Namely if
$\alpha:X\to Y$ is an equivalence of homotopy, then there exists a
bijection $S_{n}(\alpha):S_{n}(X)\to S_{n}(Y)$ defined by setting:
\begin{equation*}\label{22}
    S_{n}(\alpha)([\widetilde{\phi}_n] )=([\widetilde{\gamma_{n}^{\alpha}}\circ\widetilde{\phi}_n]
    ).
\end{equation*}
where $\widetilde{\gamma_{n}^{\alpha}}:\mathrm{Coker }\,b_{n+1}\to
\mathrm{Coker }\,b'_{n+1}$ is the quotient homomorphism induced by
$\gamma_{n}^{\alpha}:\Gamma_{n}^{X}\to \Gamma_{n}^{Y}$. Note that
the following commutative diagram:

\begin{picture}(300,90)(10,30)
\put(90,102){$\mathrm{Ext }^1_{\Bbb Z}(H_n(X), \mathrm{Coker
    }\,b_{n+1})\overset{\theta_{n}}{\twoheadrightarrow} \mathrm{Ext }(\ker b_n, \mathrm{Coker }\,b_{n+1})$}
\put(90,48){$\mathrm{Ext }^1_{\Bbb Z}(H_n(Y), \mathrm{Coker
    }\,b'_{n+1})\overset{\theta'_{n}}{\twoheadrightarrow} \mathrm{Ext }(\ker b'_n, \mathrm{Coker }\,b'_{n+1})$}
 \put(130,98){$\vector(0,-1){39}$}
 \put(132,78){\scriptsize $\cong$}
\put(277,78){\scriptsize $\cong$}
\put(275,98){$\vector(0,-1){39}$}
\end{picture}

\noindent assures that
$([\gamma_{n}^{\alpha}\circ\widetilde{\phi}_n]\in S_{n}(Y).$	

\end{remark}
We can now tackle our main theorem.

\section{The classification of CW-complexes}

\subsection{Preliminary settings}

We now go back to the problem mentioned in section 2: suppose given
two simply connected CW-complexes $X$ and $Y$ and a graded
homomorphism  $f_* : H_*(X)\rightarrow H_*(Y)$. What can we say
about the existence of a cellular map $\alpha : X\rightarrow Y$ with
$H_*(\alpha)=f_*$?

\noindent We actually need to know $f_* $ already at the chain
complexes level, say $\xi_* : C_*X \rightarrow C_*Y$ a
representative, the existence of which is certified by the homotopy
extensions theorem (see \cite{D}).

\noindent We shall proceed by induction. So we first define:
\begin{definition}
\label{d3}
The map $\alpha^n : X^n \rightarrow Y^n$ is an $n$-realization of $f_*$ if $H_{\leq n-1}(\alpha)=f_{\leq n-1}$
and $C_n\alpha^n|_{\mbox{ker }d_n} = \xi_n|_{\mbox{ker }d_n}$. We denote by $A_n=\{\alpha^n\}$ the set of all $n$-realization of $f_*$.
\end{definition}
We need further some compatibility with the Whitehead exact
sequences.
\begin{definition}
\label{d4} The pair $(f_*, \gamma_*)$ where $f_* : H_*(X)
\rightarrow H_*(Y)$ and $\gamma_* : \Gamma^X_* \rightarrow
\Gamma^Y_*$ are graded group maps  is called a morphism from
WES($X$) towards WES($Y$) if the following two property are
satisfied:

1) There exists a graded homomorphism  $\Omega_* : \pi_*(X)
\rightarrow \pi_*(Y)$ making the diagram (2) commute.

2) For every $n\geq 2$, if there exists an $n$-realization of $f_*$, then  $A_n$ contains some $\alpha^n$ with $\gamma_n= \gamma_n^{\alpha^n}$.
\end{definition}
\begin{example}
\label{e1} If $\theta : X \rightarrow Y$ is a map of CW-complexes,
it induces a morphism $(H_*(\theta), \gamma_*^{\theta})$ of
Whitehead's certain exact sequences. Obviously, this is a morphism
in the meaning of definition \ref{d4}.
\end{example}
\begin{example}
\label{e2} Theorem \ref{t1} gives a non-trivial illustration of
definition \ref{d4}. Moreover, it is an example of the forthcoming
definition \ref{d5}; this fact justifies the presentation of our
main theorem as a generalization of the result of Whitehead recalled
in subsection 2.3.
\end{example}
\begin{example}
\label{e3} In this example we are motivated by the following results
due to Anick \cite{A1,A2}.

\noindent  Let $R$ be a subring of $\Bbb Q$ and let $p$ the least
prime number which is not a unit in $R$. A free differential graded
Lie algebra is called $n$-mild if it generated by the elements with
degree $i$, where $n\leq i\leq np-1$.

\noindent Let $CW^{np}_{n}(R)/_{\simeq}$ and
$DGL^{np}_{n}/_{\simeq}$ denote the homotopy category of
$R$-localized, $n$-connected, $np$-dimensional CW-complexes and of
$n$-mild free dgl's, respectively.  Anick  has proved that the
universal enveloping algebra functor $U:DGL^{np}_{n}\to
HAH^{np}_{n}$ induces an isomorphism on the homotopy categories.
Here $HAH^{np}_{n}$ is the category of $n$-mild Hopf algebras up to
homotopy. Thus we have an equivalence of categories
$L:CW^{np}_{n}(R)/_{\simeq}\to DGL^{np}_{n}/_{\simeq}$ by composing
the Adams-Hilton model $L:CW^{np}_{n}/_{\simeq}\to
HAH^{np}_{n}/_{\simeq}$ with $U^{-1}$. Moreover if
$L(X)=(L(V),\partial )$, then for every $i<D=min(n+2p-3,np-1)$ we
have:
$$\pi _{i}(X)\otimes R\cong H_{i-1}(L(V),\partial )\text{ and }%
H_{i}(X,R)\cong H_{i}(s^{-1}V,d))\,\,\,\,\,,\,\,\, \forall i<D.$$
here $s^{-1}$ denotes the desuspension graded homomorphism.

\noindent Define $\Gamma_{i}^{L(V)}=\mathrm{Im }\big( i_n :
H_i(L(V_{i-1}) \rightarrow H_i(L(V_{i-1}))\big)$, for all $i<D$.
 Because of the above equivalence, we can say that the abelian groups
$\Gamma_{i}^{X}$ and $\Gamma_{i}^{L(V)}$ are isomorphic.

Now let  $n=1$, $p=7$ and let $X$ be an object in
$CW^{5}_{1}(R)/_{\simeq}$. By putting $H_{i}=H_{i}(X)$, for $2\leq
i\leq 5$, and by using proposition 3.4 in \cite{Ben4} we  derive
that:
\begin{equation}\label{222}
    \Gamma_{3}^{X}=[H_{2},H_{2}]\,\,\,\,\,\,,\,\,\,\,\,\,\Gamma_{4}^{X}=[H_{3},H_{2}]\oplus
    [H_{2},\mathrm{Coker}\,b_{5}]
\end{equation}
Therefore WES(X) can be written as follows:
$$H_{5}\overset{b_{5}}{\longrightarrow}[H_{3},H_{2}]\oplus [H_{2},\mathrm{Coker}\,b_{4}]\to\pi_{4}(X)\to H_{4}\overset{b_{4}}{\longrightarrow}[H_{2},H_{2}]\to\pi_{3}(X)\twoheadrightarrow H_{3}$$
Likewise  it is also shown    that if $\alpha:X\to Y$ is a morphism
in $CW^{5}_{1}(R)/_{\simeq}$, then in the following diagram:

\begin{picture}(300,90)(55,10)
\put(65,80){$H_{5}\overset{b_{5}}{\longrightarrow}[H_{3},H_{2}]\oplus
[H_{2},\mathrm{Coker}\,b_{4}]\to\pi_{4}(X)\to
H_{4}\overset{b_{4}}{\longrightarrow}[H_{2},H_{2}]\to\pi_{3}(X)\twoheadrightarrow
H_{3}$} \put(255,76){$\vector(0,-1){46}$}
\put(155,76){$\vector(0,-1){46}$}
 \put(65,20){$H'_{5}\overset{b'_{4}}{\longrightarrow}[H'_{3},H'_{2}]\oplus [H'_{2},\mathrm{Coker}\,b'_{4}]\to\pi_{4}(Y)\to H'_{4}\overset{b'_{4}}{\longrightarrow}[H'_{2},H'_{2}]\to\pi_{3}(Y)\twoheadrightarrow H'_{3}$}
\put(69,76){$\vector(0,-1){46}$} \put(71,50){\scriptsize
$H_{5}(\alpha)$} \put(307,50){\scriptsize  $H_{4}(\alpha)$}
\put(305,76){$\vector(0,-1){46}$} \put(158,50){\scriptsize
$\gamma_{4}^{\alpha}$} \put(256,50){\scriptsize $\pi_{4}(\alpha)$}
\put(358,50){\scriptsize
$\gamma_{3}^{\alpha}$}\put(356,76){$\vector(0,-1){46}$}
\put(405,76){$\vector(0,-1){46}$}\put(406,50){\scriptsize
$\pi_{3}(\alpha)$} \put(457,50){\scriptsize  $H_{3}(\alpha)$}
\put(455,76){$\vector(0,-1){46}$}
\end{picture}

\noindent where the bottom sequence is WES(Y), we have:
$$\gamma_{3}^{\alpha}=[H_{2}(\alpha),H_{2}(\alpha)]\,\,\,\,\,\,,\,\,\,\,\,\,\gamma_{4}^{\alpha}=[H_{3}(\alpha),H_{2}(\alpha)]\oplus
    \big[H_{2}(\alpha),\overline{[H_{2}(\alpha),H_{2}(\alpha)]}\big].$$
    Here
    $\overline{[H_{2}(\alpha),H_{2}(\alpha)]}:\mathrm{Coker}\,b_{4}\to
    \mathrm{Coker}\,b'_{4}$ is the quotient homomorphism induced by
    $[H_{2}(\alpha),H_{2}(\alpha)]$.

Hence if $X$ and $Y$ are in $CW^{5}_{1}(R)/_{\simeq}$, then
morphisms from WES(X) towards WES(Y) defined in \ref{d4} can be
characterized as follows:

\noindent they are homomorphisms $f_{i}:H_{i}(X)\to H_{i}(Y),
i=2,3,4,5,$ for which there exist homomorphisms $\Omega_{3},
\Omega_{4}$ making the following diagram commutes:

\begin{picture}(300,90)(55,10)
\put(65,80){$H_{5}\overset{b_{5}}{\longrightarrow}[H_{3},H_{2}]\oplus
[H_{2},\mathrm{Coker}\,b_{4}]\to\pi_{4}(X)\to
H_{4}\overset{b_{4}}{\longrightarrow}[H_{2},H_{2}]\to\pi_{3}(X)\twoheadrightarrow
H_{3}$} \put(255,76){$\vector(0,-1){46}$}
\put(155,76){$\vector(0,-1){46}$}
 \put(65,20){$H'_{5}\overset{b'_{4}}{\longrightarrow}[H'_{3},H'_{2}]\oplus [H'_{2},\mathrm{Coker}\,b'_{4}]\to\pi_{4}(Y)\to H'_{4}\overset{b'_{4}}{\longrightarrow}[H'_{2},H'_{2}]\to\pi_{3}(Y)\twoheadrightarrow H'_{3}$}
\put(69,76){$\vector(0,-1){46}$} \put(71,50){\scriptsize $f_{5}$}
\put(307,50){\scriptsize  $f_{4}$} \put(305,76){$\vector(0,-1){46}$}
\put(158,50){\scriptsize $\gamma_{4}$} \put(256,50){\scriptsize
$\Omega_{4}$} \put(358,50){\scriptsize
$\gamma_{3}$}\put(356,76){$\vector(0,-1){46}$}
\put(405,76){$\vector(0,-1){46}$}\put(406,50){\scriptsize
$\Omega_{3}$} \put(457,50){\scriptsize  $f_{3}$}
\put(455,76){$\vector(0,-1){46}$} \put(488,50){\scriptsize $(A)$}
\end{picture}

\noindent where:
$$\gamma_{3}=[f_{2},f_{2}]\,\,\,\,\,\,,\,\,\,\,\,\,\gamma_{4}=[f_{3},f_{2})]\oplus
    \big[f_{2},\overline{[f_{2},f_{2}]}\big].$$
    Here
    $\overline{[f_{2},f_{2}]}:\mathrm{Coker}\,b_{4}\to
    \mathrm{Coker}\,b'_{4}$ is the quotient homomorphism induced by
    $[f_{2},f_{2}]$.

\end{example}
Recall that we are given two CW-complexes $X$ and $Y$ and a map
$(f_*,  \gamma_*)$ between their  Whitehead's certain exact
sequences. Let us denote by:
$$ (f_n)^* : \mbox{Ext }^1_{\mathbb{Z}}(H_n(Y), \mbox{Coker }b'_{n+1})\rightarrow \mbox{Ext }^1_{\mathbb{Z}}(H_n(X), \mbox{Coker }b'_{n+1})$$
\noindent the obvious map induced by $f_n$. If $(f_*, \gamma_*)$ is a morphism, $\gamma_n : \Gamma_n^X\rightarrow  \Gamma_n^Y$ defines a quotient morphism $\widetilde{\gamma}_n : \mbox{Coker }b_{n+1}\rightarrow \mbox{Coker }b'_{n+1}$, and therefore a group morphism
$$ (\widetilde{\gamma}_n )_* : \mbox{Ext }^1_{\mathbb{Z}}(H_n(X), \mbox{Coker }b_{n+1}) \rightarrow \mbox{Ext }^1_{\mathbb{Z}}(H_n(X), \mbox{Coker }b'_{n+1})$$
\begin{definition}
\label{d5} $(f_*,  \gamma_*)$ is a strong morphism if there exist
 $[\widetilde{\phi}_n]\in S_{n}(X)$  and
$[\widetilde{\phi'}_{n}]\in  S_{n}(Y)$ such that:
\begin{equation}\label{23}
(f_n)^*([\widetilde{\phi'}_n] )= (\widetilde{\gamma}_n
)_*([\widetilde{\phi}_n])\,\,\,\,,\,\,\, \forall n\geq 2
\end{equation}
\end{definition}
\begin{remark}\label{r1}
If $X$ and $Y$ are two CW-complexes whose homology $H_*(X)$ and
$H_*(Y)$ are $\mathbb{Z}$-free, any morphism between their
respective Whitehead exact sequences is strong.
\end{remark}

For proving the main theorem,  we shall explicit condition
(\ref{23}) on representatives. As a preliminary, we achieve that in
the following subsection.
\subsection{Expliciting classes in $\mbox{Ext}$-groups}

 First let us choose the free resolution of $H_n(X)$ (resp.
$H_n(Y)$) given in (\ref{21}):
 $\mbox{Im }d_{n+1} \stackrel{\kappa_n}{\rightarrowtail}  \mbox{ker }d_n \twoheadrightarrow H_n(X)$ (resp $\mbox{Im }d'_{n+1} \stackrel{\kappa'_n}{\rightarrowtail}  \mbox{ker }d'_n \twoheadrightarrow H_n(Y)$). The cycle $\phi_n$ and its quotient $\widetilde{\phi_n}$ which represents the class $[\widetilde{\phi_n}]$ can be inserted in the following diagram (resp. for $\phi'_n$)

 \begin{picture}(300,155)(-20,-15)
\put(20,98){$\vector(0,-1){39}$} \put(10,102){$
\mathrm{Im}\,d_{n+1}\overset{\kappa_{n}}{\rightarrowtail }\ker
d_{n}\twoheadrightarrow H_{n}(X)$}
\put(210,50){$(\mathrm{Im}\,\kappa'_{n+1})'\overset{\kappa'_{n}}{\rightarrowtail
}\ker d'_{n}\twoheadrightarrow H_{n}(Y)$}
\put(210,102){$(\mathrm{Im}\,d_{n+1})'\overset{\kappa_{n}}{\rightarrowtail
}\ker d_{n}\twoheadrightarrow H_{n}(X)$}
\put(0,-01){$\mathrm{Coker}\,b'_{n+1}\underset{pr'}{\twoheadleftarrow}
\Gamma_{n}^{Y}$} \put(0,50){$\mathrm{Coker}
\,b_{n+1}\underset{pr}{\twoheadleftarrow} \Gamma_{n}^{X}$}
\put(3,78){\scriptsize $\widetilde{\phi_{n}}$}
\put(20,47){$\vector(0,-1){39}$}\put(56,78){\scriptsize $\phi_{n}$}
\put(208,-01){$\mathrm{Coker}
\,b'_{n+1}\underset{pr'}{\twoheadleftarrow} \Gamma_{n}^{Y}$}
\put(225,98){$\vector(0,-1){39}$} \put(76,47){$\vector(0,-1){38}$}
\put(225,47){$\vector(0,-1){39}$}
 \put(0,25){\scriptsize $\widetilde{\gamma_{n}}$}
 \put(78,25){\scriptsize $\gamma_{n}$}
 \put(204,25){\scriptsize $\widetilde{\phi'_{n}}$} \put(266,25){\scriptsize $\phi'_{n}$}
\put(204,78){\scriptsize $\xi_{n+1}$} \put(349,78){\scriptsize
$f_{n}$} \put(25,97){$\vector(4,-3){48}$}
\put(232,47){$\vector(4,-3){48}$} \put(345,98){$\vector(0,-1){37}$}
\put(395,50){(5)}
\end{picture}

\noindent Notice that these commutative diagrams hold for any
lifting $\phi_n$ (resp. $\phi'_n$) of $\widetilde{\phi_n}$ (resp.
$\widetilde{\phi'_n}$) - not only for the representatives defined by
formula (\ref{23}).

\noindent   Now we define $(f_n)^*$ and $(\widetilde{\gamma_n})_*$
on cycles by using the free resolution (21); taking classes again we
may write down:
 $$(\widetilde{\gamma_n})_*([\widetilde{\phi_n}])=[\widetilde{\gamma_n} \circ  \widetilde{\phi_n}] \mbox{ and } (f_n)^*( [\widetilde{\phi'_n}])=
 [ \widetilde{\phi'_n} \circ \xi_{n+1}].$$

\noindent  With these description, condition (\ref{23}) turns out to
be:
 $$[\widetilde{\gamma}_n \circ  \widetilde{\phi_n}  - \widetilde{\phi'_n} \circ \xi_{n+1}] = 0 \mbox{ in }
 \mbox{Ext }^1_{\mathbb{Z}}(H_n(X), \mbox{Coker }b'_{n+1}).$$
 Going back to cycles we deduce the existence of a group morphism $\widetilde{h_n} : \mbox{ker }d_n \rightarrow \mbox{Coker }b'_{n+1}$ satisfying:
 \begin{equation}\label{25}
\widetilde{\gamma}_n \circ  \widetilde{\phi_n}  -
\widetilde{\phi'_n} \circ \xi_{n+1} = \widetilde {h_n} \circ
\kappa_n
 \end{equation}

\noindent  As $\mbox{ker }d_n$ is free, we can find a morphism $h_n
\rightarrow \Gamma^Y_n$ which lifts $\widetilde {h}_n$:

  \begin{picture}(300,90)(-35,30)
\put(120,98){$\vector(0,-1){42}$} \put(124,98){$\vector(3,-2){60}$}
 \put(156,78){\scriptsize $\widetilde{h}_{n}$}
\put(110,102){$\ker d_{n}$}
\put(110,45){$\mathrm{Coker}\,b'_{n+1}\overset{pr}{\twoheadleftarrow}\Gamma_{n}^{Y}\subseteq
\pi_{n}(Y^{n})$} \put(107,78){\scriptsize $h_{n}$}
\end{picture}

 \noindent Thus, lifting equation (\ref{25}) - see diagrams (24 ), we get:
 \begin{equation}
 \label{26}
 \mbox{Im }(\gamma_n \circ  \phi_n  - \phi'_n \circ \xi_{n+1} - h_n \circ \kappa_n) \subset \mbox{Im
 }b'_{n+1}
 \end{equation}
Recalling the splitting map (\ref{18}):
$$t_{n+1} :  \mbox{Im }d_{n+1} \oplus \mbox{ker }d_{n+1}  \stackrel{\cong}{\rightarrow} C_{n+1}X$$
and the definition (\ref{2}) of  $\Gamma^Y_n$, we get the following
lifting $g_n$ of $h_n$:

\begin{picture}(300,165)(-45,-40)
\put(110,98){$\vector(0,-1){42}$} \put(122,97){$\vector(4,-3){129}$}
 \put(180,58){\scriptsize $g_{n}$}
\put(105,102){$C_{n}X$}
\put(55,-17){$\pi_{n}(Y^{n})\supseteq\Gamma_{n}^{Y}$}
\put(230,-20){$\vector(-1,0){100}$}
\put(237,-20){$\pi_{n}(Y^{n-1})$} \put(83,78){\scriptsize
$(t_{n})^{-1}$} \put(73,45){$\mbox{Im }d_n \oplus  \mbox{ker }d_n$}
\put(110,40){$\vector(0,-1){42}$} \put(94,20){\scriptsize $h_{n}$}
\put(180,-28){\scriptsize $i_{n}$}
\end{picture}

\noindent Finally formula (\ref{26}) becomes
\begin{equation}\label{27}
\mbox{Im }(\gamma_n \circ  \phi_n  - \phi'_n \circ \xi_{n+1} - i_n
\circ g_n \circ d_{n+1}) \subset \mbox{Im }b'_{n+1}
\end{equation}
where we use the fact that $\kappa_n \circ (t_{n+1})^{-1}= d_{n+1}$.
Once again we emphasize that $\phi_n$ (resp. $\phi'_n$) is any
lifting of $\widetilde{\phi_n}$ (resp. $\widetilde{\phi'_n}$).
\subsection{The main theorem}
We can now formulate and prove the following:
\begin{theorem}
\label{t3} Let $X$ and $Y$ be two simply connected CW complexes and
$(f_*, \gamma_*) : \mbox{WES }(X) \rightarrow  \mbox{WES }(Y)$ a
strong homomorphism. Then there exists a cellular map $\alpha : X
\rightarrow Y$ such that $H_*(\alpha) = f_*$.
\end{theorem}
As an immediate consequence of the Whitehead theorem, we get:
\begin{corollary}
\label{c1} Let $X$ and $Y$ be two simply connected CW complexes. If
$\mbox{WES }(X)$ and $ \mbox{WES }(Y)$ are strongly isomorphic, then
$X$ and $Y$ are homotopic.
\end{corollary}
\begin{remark}
\label{r3} If $H_*(X)$ and $H_*(Y)$ are free abelian groups, the
assertion "strong" in the above corollary is automatically
satisfied.
\end{remark}
The last pages are now devoted to the proof of theorem \ref{t3}.
\begin{remark}
\label{r4} For the proof of theorem \ref{t3}  we need  the following
 elementary fact:

Let $X$ and $Y$ be CW-complexes and $F : X^n \rightarrow Y^n$ a map
between their $n$-skeleta. If $\xi_{n+1}: C_{n+1}(X) \rightarrow
C_{n+1}(Y)$ is a morphism such that the following diagram commutes:

\begin{picture}(300,90)(-30,30)
\put(72,100){$ C_{n+1}X
\hspace{1mm}\vector(1,0){93}\hspace{1mm}C_{n+1}Y$}
 \put(89,76){$\beta_{n+1}$} \put(214,76){$\beta'_{n+1}$}
\put(87,97){$\vector(0,-1){38}$} \put(212,96){$\vector(0,-1){38}$}
\put(145,103){ $\rho_{n+1}$} \put(140,52){ $\pi_{n}(F)$}
\put(80,48){$ \pi_{n}(X^{n})\hspace{1mm}\vector(1,0){80}\hspace{1mm}
\pi_{n}(Y^{n})$}
\end{picture}

 then $F$ can be
extended in $G: X^{n+1} \rightarrow Y^{n+1}$ with $C_{n+1}(G) =
\rho_{n+1}$.
\end{remark}
\begin{proof}(of theorem \ref{d3})\\
The proof goes recursively. At each step, we wish to set on the data
in such a way that we can apply remark \ref{r4}.

Let $X$ and $Y$ be two simply connected CW complexes. By hypothesis we are given the following commutative diagram

\begin{picture}(300,90)(25,10)
\put(65,80){$\cdots \rightarrow H_{n+1}(X,\Bbb
Z)\overset{b_{n+1}}{\longrightarrow }\Gamma^{X} _{n}\longrightarrow
\pi_{n}(X)\longrightarrow H_{n}(X,\Bbb
Z)\overset{b_{n}}{\longrightarrow }\cdots$}
\put(235,76){$\vector(0,-1){46}$} \put(186,76){$\vector(0,-1){46}$}
 \put(65,20){$\cdots \rightarrow
H_{n+1}(Y,\Bbb Z)\overset{b'_{n+1}}{\longrightarrow }\Gamma^{Y}
_{n}\longrightarrow \pi_{n}(Y)\longrightarrow H_{n}(Y,\Bbb
Z)\overset{b'_{n}}{\longrightarrow }\cdots$}
\put(110,76){$\vector(0,-1){46}$} \put(111,50){\scriptsize
$f_{n+1}$} \put(291,50){\scriptsize  $f_{n}$}
\put(290,76){$\vector(0,-1){46}$} \put(188,50){\scriptsize
$\gamma_{n}$} \put(236,50){\scriptsize $\Omega_{n}$}
\end{picture}

\noindent and a chain map $\xi_* : C_*X \rightarrow C_*Y$ whose
homology is $H_*(\xi)=f_*$.

\noindent   Suppose now we have already constructed an
$n$-realization of $(f_*,  \gamma_*)$. By definition \ref{d3} it
means there exists a map $\alpha^n : X^n \rightarrow Y^n$  whith the
following properties:
\begin{equation}\label{28}
H_{\leq n-1}(\alpha^n)=f_{\leq n-1}, C_n\alpha^n|_{\mbox{ker
}d_n}=\xi|_{\mbox{ker }d_n} \mbox{and } \gamma_n=\gamma_n^{\alpha^n}
\end{equation}
Now let us consider the following diagram:

\begin{picture}(300,190)(-30,-65)
\put(82,100){$ C_{n+1}X
\hspace{1mm}\vector(1,0){120}\hspace{1mm}C_{n+1}Y$}
 \put(89,76){\scriptsize $\beta_{n+1}$} \put(250,76){\scriptsize $\beta'_{n+1}$}
\put(87,97){$\vector(0,-1){36}$} \put(248,96){$\vector(0,-1){36}$}
\put(165,103){\scriptsize $\xi_{n+1}$} \put(160,52){\scriptsize
$\scriptsize \pi_{n}(\alpha^{n})$} \put(60,48){$
\pi_{n}(X^{n})\hspace{1mm}\vector(1,0){140}\hspace{1mm}\pi_{n}(Y^{n})$}
\put(42,-4){$ \Gamma_{n}^{X}\oplus \ker
\beta_{n}\hspace{1mm}\vector(1,0){135}\hspace{1mm}
\Gamma_{n}^{Y}\oplus \ker \beta'_{n}$} \put(89,24){\scriptsize
$\mu_{n}$} \put(250,24){\scriptsize $\mu'_{n}$}
\put(75,24){\scriptsize $\cong$} \put(238,24){\scriptsize $\cong$}
\put(87,45){$\vector(0,-1){36}$} \put(248,45){$\vector(0,-1){36}$}
\put(130,-15){\scriptsize $\gamma^{\alpha^{n}}_{n}\oplus
C_{n}\alpha^{n}\mid_{\ker \beta'_{n}}$}
 \put(-17,-50){$ C_{n}X
\hspace{1mm}\vector(1,0){330}\hspace{1mm}C_{n}Y$}
\put(78,42){$\vector(-1,-1){80}$} \put(80,95){$\vector(-2,-3){88}$}
\put(15,35){\scriptsize $d_{n+1}$}\put(38,16){\scriptsize $j_{n}$}
\put(265,42){$\vector(1,-1){80}$}
\put(268,95){$\vector(2,-3){88}$}\put(165,-47){\scriptsize
$\xi_{n}$} \put(312,35){\scriptsize
$d'_{n+1}$}\put(294,16){\scriptsize $j'_{n}$} \put(385,50){(6)}
\end{picture}

\noindent In this diagram,  both squares do not commute, but both
side triangles and both trapeza do. (These various commutative
subdiagrams translate the respective definitions of the differential
$d_*$ of the chain complex $C_*X$, of the  chain map $\xi_* : C_*X
\rightarrow C_*Y$ and  of the cellular map $\alpha^n :
X^n\rightarrow Y^n$. The lower (non-commutative) square is merely
diagram (4)).

\noindent  We shall show that we can pertubate both $\alpha^n$ and
$\xi_{n+1}$ so that we can apply remark \ref{d5} and still satisfy
the recursive hypothesis, proving the "$n$ to $n+1$ step". Let us
now explain the details.

\noindent  To begin with let us examine the big square in diagram
(6); more precisely, we calculate
$$(\gamma^{\alpha^n}_n \oplus C_n\alpha^n|_{\mbox{ker }\beta_n}) \circ \mu_n \circ \beta_{n+1} -  \mu'_n \circ \beta'_{n+1} \circ \xi_{n+1}$$
Let us make use of the matrix setting given in (\ref{1111}). In this
setting the first summand has the following expression:

\begin{equation}\label{30}
\left(
\begin{array}{cc}
\gamma^{\alpha^n}_n & 0 \\
0 & C_n\alpha^n|_{\mbox{ker }\beta_n} \\
\end{array}
\right)
 \circ
\left(
\begin{array}{cc}
\phi_n & b_{n+1}\circ pr_{n+1} \\
d_{n+1}\circ t_{n+1}|_{\mbox{Im }d_{n+1}} & 0 \\
\end{array}
\right)
$$
$$=
\left(
\begin{array}{cc}
\gamma^{\alpha^n}_n \circ \phi_n & \gamma^{\alpha^n}_n \circ b_{n+1}\circ pr_{n+1} \\
C_n\alpha^n|_{\mbox{ker }\beta_n} \circ d_{n+1}\circ t_{n+1}|_{\mbox{Im }d_{n+1}} & 0 \\
\end{array}
\right)
\end{equation}

For the second summand:

\begin{equation}\label{31}
\left(
\begin{array}{cc}
\phi'_n & b'_{n+1}\circ pr'_{n+1} \\
d'_{n+1}\circ t'_{n+1}|_{\mbox{Im }d'_{n+1}} & 0 \\
\end{array}
\right)
 \circ
 \left(
\begin{array}{cc}
pr_1 \circ (t'_{n+1})^{-1} \circ \xi_{n+1} \circ t_{n+1}|_{\mbox{Im }d_{n+1}} & 0 \\
pr_2 \circ (t'_{n+1})^{-1} \circ \xi_{n+1} \circ t_{n+1}|_{\mbox{Im }d_{n+1}} & \xi_{n+1}|_{\mbox{ker  }d_{n+1}} \\
\end{array}
\right)
$$
$$=
 \left(
\begin{array}{cc}
\phi'_n \circ pr_1 \circ (t'_{n+1})^{-1} \circ \xi_{n+1} \circ t_{n+1}|_{\mbox{Im }d_{n+1}} + \Delta_{n+1}  &  b_{n+1}\circ pr'_{n+1} \circ \xi_{n+1}|_{\mbox{ker  }d_{n+1}}\\
d'_{n+1}\circ t'_{n+1}|_{\mbox{Im }d'_{n+1}} \circ pr_1 \circ (t'_{n+1})^{-1} \circ \xi_{n+1} \circ t_{n+1}|_{\mbox{Im }d_{n+1}} & 0 \\
\end{array}
\right)
\end{equation}
where $\Delta_{n+1}  = b_{n+1}\circ pr'_{n+1} \circ pr_2 \circ
(t'_{n+1})^{-1} \circ \xi_{n+1} \circ t_{n+1}|_{\mbox{Im }d_{n+1}}$.

\begin{remark}\label{r6}
As we noticed in paragraph 5.2 both maps $\phi'_n \circ pr_1 \circ
(t'_{n+1})^{-1} \circ \xi_{n+1} \circ t_{n+1}|_{\mbox{Im }d_{n+1}} +
\Delta_{n+1} \mbox{and }  \phi'_n \circ pr_1 \circ (t'_{n+1})^{-1}
\circ \xi_{n+1} \circ t_{n+1}|_{\mbox{Im }d_{n+1}} : \mbox{Im
}d_{n+1} \rightarrow \Gamma^Y_n$ project identically onto
$\mbox{Coker }b'_{n+1}$. So the explicit formula (\ref{27}) holds -
just choose an other adequate $g_n$:
\begin{equation}\label{32}
\mbox{Im }(\gamma_n \circ  \phi_n  - (\phi'_n \circ pr_1 \circ
(t'_{n+1})^{-1} \circ \xi_{n+1} \circ t_{n+1}|_{\mbox{Im }d_{n+1}}+
\Delta_{n+1} )- i_n \circ g_n \circ d_{n+1}) \subset \mbox{Im
}b'_{n+1}.
\end{equation}
\end{remark}
 Determining the lack of commutativity of the big square
in diagram (6), we have to calculate the difference of matrix
(\ref{30}) and matrix (\ref{31}); we obtain:

\begin{equation}\label{33}
\left(
\begin{array}{cc}
\gamma_n \circ  \phi_n  - (\phi'_n \circ pr_1 \circ (t'_{n+1})^{-1} \circ \xi_{n+1} \circ t_{n+1}|_{\mbox{Im }d_{n+1}}  + \Delta_{n+1} ) & 0 \\
0 & 0 \\
\end{array}
\right)
\end{equation}
We used here:

\noindent For the $0$ on the first line, the hypothesis that
$\gamma^{\alpha^n}_n = \gamma_n$; For the $0$ on the second line the
hypothesis that $C_n\alpha^n|_{\mbox{ker }\beta_n} =
C_n\alpha^n|_{\mbox{ker }d_n} = \xi_n|_{\mbox{ker }d_n}$.

\noindent We resume our calculation:
\begin{eqnarray}
&&(\gamma^{\alpha^n}_n \oplus C_n\alpha^n|_{\mbox{ker }\beta_n}) \circ \mu_n \circ \beta_{n+1} -  \mu'_n \circ \beta'_{n+1} \circ \xi_{n+1} \nonumber\\
&=& \gamma_n \circ  \phi_n  - (\phi'_n \circ pr_1 \circ
(t'_{n+1})^{-1} \circ \xi_{n+1} \circ t_{n+1}|_{\mbox{Im }d_{n+1}} +
\Delta_{n+1} )\label{34}
\end{eqnarray}
 Now we focus on the lack of commutativity of the lower square and
compose equation (\ref{17})) on the right by $\beta_{n+1}$:
\begin{eqnarray}
[\mu'_n \circ  \pi_n(\alpha^n) - (\gamma^{\alpha^n}_n \oplus
C_n\alpha^n|_{\mbox{ker }\beta_n} )\circ \mu_n ] \circ  \beta_{n+1}
\hspace{-3mm}&= &\hspace{-3mm}[\sigma'_n \circ j'_n \circ \pi_n(\alpha^n)   -  \pi_n(\alpha^n) \circ \sigma_n \circ j_n] \circ \beta_{n+1}\nonumber\\
&=&\hspace{-3mm}(\sigma'_n \circ \xi_n     -  \pi_n(\alpha^n) \circ \sigma_n) \circ j_n \circ \beta_{n+1}\nonumber\\
&= &\hspace{-3mm}( \sigma'_n \circ \xi_n  -  \pi_n(\alpha^n) \circ
\sigma_n) \circ d_{n+1} \label{35}
\end{eqnarray}
where we used the commutations $j'_n \circ \pi_n(\alpha^n) = \xi_n
\circ j_n$ and $j_n \circ \beta_{n+1} = d_{n+1}$ we pointed to in
diagram (6).

\noindent As $C_{n+1}X$ is free we can lift the morphism  $\sigma'_n
\circ \xi_n   -   \pi_n(\alpha^n) \circ \sigma_n : C_{n+1}X
\rightarrow \Gamma^Y_n$ to $\pi_n(Y^{n-1})$ as pictured in the
following diagram:

\begin{picture}(300,95)(-45,30)
\put(110,98){$\vector(0,-1){42}$} \put(122,97){$\vector(3,-1){125}$}
 \put(188,78){\scriptsize $k_{n}$}
\put(105,102){$C_{n+1}X$}
\put(55,45){$\pi_{n}(Y^{n})\supseteq\Gamma_{n}^{Y}$}
\put(230,45){$\vector(-1,0){100}$} \put(237,45){$\pi_{n}(Y^{n-1})$}
\put(27,78){\scriptsize
$\sigma'_{n}\xi_{n}-\pi_n(\alpha^n)\circ\sigma_{n}$}
\put(180,37){\scriptsize $i_{n}$}
\end{picture}

\noindent Summing formulas (\ref{34}) and (\ref{35}) we obtain:
\begin{eqnarray}
\mu'_n \circ  \pi_n(\alpha^n) \circ  \beta_{n+1} - \mu'_n \circ
\beta'_{n+1} \circ \xi_{n+1}\hspace{-2mm}& =& \hspace{-2mm} \gamma_n
\circ \phi_n -(\phi'_n \circ pr_1 \circ (t'_{n+1})^{-1} \circ
\xi_{n+1} \circ
t_{n+1}|_{\mbox{Im }d_{n+1}}\hspace{-2mm}+ \Delta_{n+1} )\nonumber\\
& + & i_n \circ k_n \circ d_{n+1}.\label{36}
\end{eqnarray}
As $\mu'_n$ is an isomorphism we can rewrite this equation
\begin{eqnarray}
  \pi_n(\alpha^n) \circ  \beta_{n+1} - \beta'_{n+1} \circ \xi_{n+1}
\hspace{-2mm}&=& \hspace{-2mm} (\mu'_n)^{-1} \circ [\gamma_n \circ  \phi_n  - (\phi'_n \circ pr_1 \circ (t'_{n+1})^{-1} \circ \xi_{n+1} \circ t_{n+1}|_{\mbox{Im }d_{n+1}}
\hspace{-2mm}+ \Delta_{n+1} ) \nonumber\\
&+& i_n \circ k_n \circ d_{n+1}] \label{37}
\end{eqnarray}
We know first that $\gamma_n \circ  \phi_n  - (\phi'_n \circ
\xi_{n+1} + \Delta_{n+1} ) + i_n \circ k_n \circ d_{n+1}$ maps into
$\Gamma^Y_n$, second that $(\mu'_n)^{-1} $ is the identity on
$\Gamma^Y_n$. We then have :
\begin{eqnarray}
\pi_n(\alpha^n) \circ  \beta_{n+1} -  \beta'_{n+1} \circ \xi_{n+1}&=&  \gamma_n \circ  \phi_n  - (\phi'_n \circ pr_1 \circ (t'_{n+1})^{-1} \circ \xi_{n+1} \circ t_{n+1}|_{\mbox{Im }d_{n+1}} + \Delta_{n+1} ) \nonumber\\
&+& i_n \circ  k_n \circ d_{n+1}\label{38}
\end{eqnarray}
Let us now substitute this last equation in formula (\ref{32}); we
get:
\begin{equation}\label{39}
\mbox{Im }(\pi_n(\alpha^n) \circ  \beta_{n+1} - \beta'_{n+1} \circ
\xi_{n+1} - i_n \circ k_n \circ d_{n+1} - i_n \circ g_n \circ
d_{n+1}) \subset \mbox{Im }b'_{n+1}.
\end{equation}
Recalling again $j_n \circ \beta_{n+1} = d_{n+1}$ this becomes:
 \begin{equation}\label{40}
\mbox{Im }([(\pi_n(\alpha^n) - i_n \circ (k_n - g_n) \circ j_n]
\circ  \beta_{n+1} - \beta'_{n+1} \circ \xi_{n+1}) \subset \mbox{Im
}b'_{n+1}.
 \end{equation}
So there is a well defined morphism:

\begin{picture}(300,40)(-10,13)
\put(30,28){$(\mathrm{Im}\,d_{n+1})'\hspace{2mm}\vector(1,0){220}\hspace{2mm}\mathrm{Im}\hspace{2mm}b'_{n+1}
\subset \Gamma _{n}^{Y}$} \put(90,35){\scriptsize $[(\pi_n(\alpha^n)
- i_n \circ (k_n - g_n)) \circ j_n] \circ \beta_{n+1} - \beta'_{n+1}
\circ \xi_{n+1} $}
\end{picture}

\noindent  It admits a lifting $\lambda_{n+1}$ as pictured in the
diagram:

\begin{picture}(300,120)(-10,-15)
\put(320,83){$\vector(0,-1){45}$} \put(306,85){$\ker d'_{n+1}$}
\put(86,30){$\vector(4,1){216}$}
\put(30,28){$(\mathrm{Im}\,d_{n+1})'\hspace{2mm}\vector(1,0){220}\hspace{2mm}\mathrm{Im}\hspace{2mm}b'_{n+1}
\subset \Gamma _{n}^{Y}$} \put(90,15){\scriptsize $[(\pi_n(\alpha^n)
- i_n \circ (k_n - g_n)) \circ j_n] \circ \beta_{n+1} - \beta'_{n+1}
\circ \xi_{n+1} $} \put(141,63){\scriptsize $\lambda _{n+1}$}
\put(322,60){\scriptsize $\beta' _{n+1}$}
\put(320,39){$\vector(0,-1){4}$}
\end{picture}

\noindent We are now ready to define maps $\psi_n : X^n \rightarrow
Y^n$ and $\rho_{n+1} : C_{n+1}X\rightarrow C_{n+1}Y$ fitting the
conditions of remark \ref{d5}.

 \noindent First we easily define  $\rho_{n+1} = \xi_{n+1} + \lambda_{n+1}$.  As to $\psi_n$ and to begin with, we choose a map $\omega  : \bigvee_{k\in K_n} S^n_k \rightarrow Y^{n-1}$ (where $K_n$ indices the $n$-cells of $X^n$) whose homotopy class corresponds to $g_n - k_n$ via the isomorphism $\mbox{Hom }(C_nX, \pi_n(Y^{n-1})) \cong [\bigvee_{k\in K_n} S^n_k, Y^{n-1}]$.

 \noindent Let us recall that attaching cells by a map such as $f : \bigvee_{k\in K} S^n_k \rightarrow X^{n-1}$ induces the action:
 $$[\bigvee_{k\in K} S^n_k, Y^n] \times [X^n, Y^n] \stackrel {\vee}{\rightarrow} [X^n, Y^n].$$
So we thus define $\psi_n = (i_n \circ \omega) \vee \alpha^n$ and check that :
\begin{equation}\label{42}
\pi_n(\psi_n) = \pi_n(\alpha^n) - i_n \circ (k_n - g_n) \circ j_n.
\end{equation}
Finally we remark the following fact: as $\mbox{Im }\omega \subset
Y^{n-1}$ the chain morphisms induced by $\psi_n$ and $\alpha^n :
C_nX \rightarrow C_nY$ agree and therefore:
\begin{equation}\label{43}
H_{\leq n}(\psi^n) = H_{\leq n}(\alpha^n).
\end{equation}
By the very definitions of $\psi_n$, $\lambda_{n+1}$ and $\rho_{n+1}$, the following diagram commutes:

\begin{picture}(300,90)(-30,30)
\put(72,100){$ C_{n+1}X
\hspace{1mm}\vector(1,0){93}\hspace{1mm}C_{n+1}Y$}
 \put(89,76){\scriptsize  $\beta_{n+1}$} \put(214,76){\scriptsize  $\beta'_{n+1}$}
\put(87,97){$\vector(0,-1){38}$} \put(212,96){$\vector(0,-1){38}$}
\put(145,103){\scriptsize  $\xi_{n+1}$} \put(140,52){\scriptsize
$\pi_{n}(\psi _{n})$} \put(80,48){$
\pi_{n}(X^{n})\hspace{1mm}\vector(1,0){80}\hspace{1mm}
\pi_{n}(Y^{n})$}
\end{picture}

\noindent  So, by remark \ref{d5}, there exists a map $\alpha^{n+1}:
X^{n+1} \rightarrow Y^{n+1} $ which extends $\psi_n$: $H_{\leq
n}(\alpha^{n+1}) = H_{\leq n}(\psi_n)$) and $H_{n+1} (\alpha^{n+1})
= \rho_{n+1}|_{\mbox{ker }d_{n+1}}$.

\noindent  Now these two equations become: First $H_{\leq
n}(\alpha^{n+1}) = H_{\leq n}(\alpha^{n})$ by equation (\ref{43}).
Second $H_{n+1} (\alpha^{n+1}) = C_{n+1}\alpha^{n+1}|_{\mbox{ker
}d_{n+1}} = \xi_{n+1}|_{\mbox{ker }d_{n+1}}$, again by the very
definition of $\rho_{n+1}$ and $\lambda_{n+1}$.

\noindent So we can phrase: $\alpha_{n+1}$ is an $(n+1)$-realization
of $(f_*, \gamma_*)$, and we got the $(n+1)$th-step of our recursive
proof.
\end{proof}
\begin{example}\label{e4} Classification of simply connected,
5-dimensional and $R$-localized CW-complexes.

 \noindent Let us consider
again the example \ref{e3}. Let $X,Y\in CW^{5}_{1}(R)/_{\simeq}$,
due to corollary \ref{c1} we can say that  if $f_{i}:H_{i}(X)\to
H_{i}(Y), i=2,3,4,5,$  are  homomorphisms for which there exist
$\Omega_{3}, \Omega_{4}$ making the
 diagram (A) commutes and  if there exist two characteristic $4$-extensions
 $[\widetilde{\phi}_4]\in S_{4}(X)\subset\mathrm{Ext }^1_{\mathbb{Z}}(H_4(X), \mathrm{Coker }b_{5})$ and
 $[\widetilde{\phi}_n]\in S_{4}(Y)\subset\mathrm{Ext }^1_{\mathbb{Z}}(H_4(Y), \mathrm{Coker
 }b'_{5})$ satisfying the equation (\ref{23}), for $n=4$. Then  $X$
 and $Y$ are homotopic.
\end{example}

\bibliographystyle{amsalpha}

\end{document}